\documentclass{amsart}

\usepackage{amssymb, latexsym}

\newtheorem{theorem}{Theorem}         
\newtheorem{corollary}{Corollary}

\newtheorem*{maintheorem}{Summary of Theorem 2}
\newtheorem*{corollary_intro}{Corollary}

\providecommand{\R}{\ensuremath{\mathbb{R}}}
\providecommand{\Z}{\ensuremath{\mathbb{Z}}}
\providecommand{\F}{\ensuremath{\mathbb{F}}}
\providecommand{\X}{\ensuremath{\mathbb{X}}}
\providecommand{\Hyp}{\ensuremath{\mathbb{H}}}
\providecommand{\G}{\Gamma}
\providecommand{\cG}{\mathcal{G}}

\providecommand{\Aut}{\mathop{\rm Aut}\nolimits}
\providecommand{\Vol}{\mathop{\rm Vol}\nolimits}
\providecommand{\Fix}{\mathop{\rm Fix}\nolimits}
\providecommand{\Lk}{\mathop{\rm Lk}\nolimits}
\providecommand{\Fix}{\mathop{\rm Fix}\nolimits}
\providecommand{\bs}{\backslash}
\providecommand{\quot}{\bs \! \bs}

\begin{document}

\title[Covolumes of uniform lattices]
 {Covolumes of uniform lattices acting on polyhedral complexes} 

\author{Anne Thomas}
\address{Department of Mathematics \\ University of Chicago \\ Chicago, Illinois 60637}
\email{athomas@math.uchicago.edu}

\begin{abstract}
Let $X$ be a polyhedral complex with finitely many isometry
classes of links. We establish a restriction on the covolumes of
uniform lattices acting on $X$. When $X$ is two-dimensional and
has all links isometric to either a complete bipartite graph or
the building for a Chevalley group of rank 2 over a field of prime
order, we obtain further restrictions on covolumes.
\end{abstract}

\maketitle

\section{Introduction}
\label{s:intro}

\noindent Let $G$ be a locally compact topological group, with suitably normalised Haar
measure $\mu$.  Let $\G \leq G$ be a \emph{uniform lattice}, that is, a
cocompact discrete subgroup. A basic question is: which covolumes $\mu(\G\bs
G)$ can occur?

When $G$ is a Lie group, either real or $p$-adic, this is a well-studied
classical problem.  See, for example,~\cite{s1:rdg},~\cite{p1:v}
and~\cite{lw1:lmc}.  A non-classical case is that of the automorphism group of
a locally finite tree, which is naturally a locally compact topological group;
see~\cite{bl1:tl}. Levich and Rosenberg~(\cite{ros1:tctl},
Chapter 9) have completely classified the covolumes of uniform lattices acting
on regular and biregular trees.

In this note, we study the higher-dimensional case of covolumes of uniform
lattices acting on locally finite $n$-dimensional polyhedral complexes, for $n
\geq 2$. Some of these complexes
are generalisations of buildings. Their automorphism groups are locally compact
when endowed with a natural topology, and the corresponding Haar measure may be
suitably normalised. See Section~\ref{s:background} for precise definitions and
this normalisation.

The main result obtained is Theorem~\ref{t:covols} below.  This theorem applies
to a large class of polyhedral complexes, those where the quotient by the full
automorphism group is finite, so that there are finitely many isometry
classes of links (the link of a vertex of a polyhedral complex is defined in
Section~\ref{ss:poly}). If $X$ is such a complex, and $\G \leq G = \Aut(X)$ is
a uniform lattice, Theorem~\ref{t:covols} may be summarised as:

\begin{maintheorem} There is a restriction on $\mu(\G\bs G)$ based only on the prime
divisors of the orders of the automorphism groups of the finite
simplicial complexes which are the links of $X$.\end{maintheorem}

The statement and proof of Theorem~\ref{t:covols} are found in
Section~\ref{ss:mainthm}, and in Section~\ref{ss:examples} we
apply Theorem~\ref{t:covols} to many examples, including some in dimensions
greater than 2. In
Section~\ref{s:onelink}, we use an analogue of
Theorem~\ref{t:covols} to establish even stronger restrictions on
covolumes of uniform lattices which act on
certain 2-dimensional complexes with only \emph{one} isometry class of
link. For example, we apply our general result to ``Bourdon
buildings'' (see Bourdon~\cite{b1:ih} and
Bourdon--Pajot~\cite{bp2:pi},~\cite{bp1:rqi}) to obtain:

\begin{corollary_intro}\label{c:intro} Let $L = K_{m,n}$ be a complete bipartite
graph, with $m \geq n \geq 2$, and let $r \geq 5$, with $r$ even
if $m \not = n$.  Let $P$ be a regular right-angled hyperbolic
$r$-gon.  Suppose $X$ is the (unique) polygonal complex with all
links isometric to $L$ and all 2-cells isometric to $P$. Then if
$\G$ is a uniform lattice in $G=\Aut(X)$,\[ \mu(\G \bs G) =
\frac{r}{mn}\left( \frac{a}{b} \right) \] where $a/b$ is rational
(in lowest terms), and
 the prime divisors of $b$ are strictly less than $m$.
\end{corollary_intro}

In fact, this result is sharp, in the sense that, as we show in~\cite{at1:cb},
every possible covolume satisfying the restrictions of this corollary actually
occurs for some uniform lattice $\G \leq G$.

\vspace{3mm} 
I would like to thank Benson Farb for suggesting
this area of research, and for his guidance and valuable advice. I
would also like to thank G. Christopher Hruska for numerous
rewarding conversations, and F. Buekenhout for his help with the 120-cell.

\section{Background}\label{s:background}

In this section we give the basic definitions of polyhedral
complexes, their automorphism groups and lattices, and describe
several methods for constructing these complexes and lattices.

\subsection{Polyhedral complexes}\label{ss:poly}

Let $\X^n$ be $S^n$, $\R^n$ or $\Hyp^n$, endowed with metrics of constant
curvature 1, 0 and $-1$ respectively. A \emph{polyhedral complex} $X$ is a
CW-complex such that: \begin{enumerate} \item each open cell of
dimension $n$ is isometric to the interior of a compact convex polyhedron in
$\X^n$; and \item for each cell $\sigma$ of $X$, the restriction of the
attaching map to each open codimension 1 face of $\sigma$ is an isometry onto
an open cell of $X$. \end{enumerate} A polyhedral complex is said to be
(piecewise) \emph{spherical}, \emph{Euclidean} or \emph{hyperbolic} if $\X^n$
is $S^n$, $\R^n$ or $\Hyp^n$ respectively.  A 2-dimensional polyhedral complex
is called a polygonal complex.

Let $x$ be a vertex of a $d$-dimensional polyhedral complex $X$. The
\emph{link} of $x$, written $\Lk(x,X)$, is the spherical simplicial complex
obtained by intersecting $X$ with a $d$-sphere of sufficiently small radius
centred at the vertex $x$.  For example, if $X$ has dimension 2, then
$\Lk(x,X)$ may be identified with the graph having vertices the 1-cells of $X$
containing $x$ and edges the 2-cells of $X$ containing $x$; an edge joins two
vertices in the link if the corresponding 2-cell contains those two 1-cells. By
rescaling so that for each $x$ the $d$-sphere around $x$ has radius say 1, we
induce a metric on each link, and we may then speak of isometry classes of
links of $X$.

There are several constructions of polyhedral complexes with specified links.

\begin{enumerate} \item Buildings of dimension $n$ have cells which are
Euclidean simplexes, and the links of all vertices isometric to a spherical
building of dimension $(n-1)$.  \item The Davis--Moussong complex~\cite{d1:bcat0}
is constructed from a metric flag complex $L$.  There exists a CAT(0)
polyhedral complex $X_L$ such that the link of every vertex in $X_L$ is
isometric to $L$. \item Ballmann--Brin~\cite{bb1:pc}, by attaching edges and
polygons step-by-step, have constructed polygonal complexes with all cells
regular Euclidean $k$-gons, $k \geq 6$, and all links the 1-skeleton of either
the $n$-simplex, the $n$-cube or a Platonic solid. \item Polyhedral complexes
may be realised as universal covers of complexes of groups (see~\cite{s1:tg},
\cite{c1:cg}, \cite{h1:cg} and \cite{bh1:ms} for the theory of complexes of
groups).  In this way, Haglund~\cite{h2:pg} has obtained similar results
to~\cite{bb1:pc}, while Benakli has constructed polygonal complexes so that
every link is isometric to one of finitely many given graphs~\cite{ben1:pc1}.
Bourdon has constructed so-called \emph{Fuchsian buildings}, which are
hyperbolic polygonal complexes, using complexes of groups~\cite{b2:if}.
\item\label{e:hp} Haglund--Paulin~\cite{hp1:ca} have constructed 3-dimensional
polyhedral complexes with each 3-cell a hyperbolic polytope and the link of
each vertex the flag complex of projective 3-space over a
finite field.  They used decompositions of buildings, which reflect
decompositions of their Coxeter systems.  \end{enumerate}

\subsection{Lattices and covolumes}\label{ss:lattcovol}

Let $G$ be a locally compact topological group with left-invariant Haar measure
$\mu$.  A discrete subgroup $\Gamma \leq G$ is called a \emph{lattice} if
the covolume $\mu(\Gamma\backslash G)$ is finite. A lattice $\G$ is
\emph{uniform} if $\G \bs G$ is compact. Let
$S$ be a left $G$-set such that for every $s \in S$, the stabiliser $G_s$ is
compact and open. Then if $\G \leq G$ is discrete, the stabilisers $\Gamma_s$
are finite. We then define the \emph{$S$-covolume} of $\Gamma$ by
\[\Vol(\Gamma \quot S)  =  \sum_{s \in \Gamma
\bs S } \frac{1}{|\Gamma_s|} \leq \infty \] The
following theorem shows that Haar measure may be normalised so
that $\mu(\G \bs G)$ equals the $S$-covolume.

\begin{theorem}[(\cite{bl1:tl}, Chapter 1)]\label{t:Scovolumes} Let $G$ be a
locally compact topological group acting on a set $S$ with compact open
stabilisers and a finite quotient $G\backslash S$.  Suppose further that $G$
admits at least one lattice.  Then there is a normalisation of the Haar
measure $\mu$, depending only on the choice of $G$-set $S$, such that for each
discrete subgroup $\Gamma$ of $G$ we have $ \mu(\G \bs G)=\Vol(\Gamma \quot S).
$ \end{theorem}

Let $X$ be a connected, locally finite, $d$-dimensional polyhedral
complex, with vertex set $V(X)$. We write $\Aut(X)$ for the group
of polyhedral isometries of $X$.  A subgroup of $\Aut(X)$ is said
to act \emph{without inversions} if its elements fix pointwise
each cell that they preserve. Let $\sigma$ be a cell of $X$.  For
$n \geq 0$, we define \emph{combinatorial balls} $B(\sigma,n)$
centred at $\sigma$ by induction. The ball $B(\sigma,0)$ is just
the cell $\sigma$, while for $n \geq 1$, $B(\sigma, n)$ is the
union of the $d$-cells of $X$ which meet $B(\sigma,n-1)$.

The group $G = \Aut(X)$ naturally has the structure of a locally
compact topological group, with a neighbourhood basis of the
identity consisting of automorphisms fixing larger and larger
combinatorial balls.  A subgroup $\Gamma \leq G$ is discrete in
this topology if and only if the stabiliser $\G_x$ is finite for
each $x \in V(X)$. By the same arguments as for tree
lattices~(\cite{bl1:tl}, Chapter 1), it can be shown that if $G\bs
X$ is finite, then a discrete subgroup $\G \leq G$ is a uniform
lattice if and
only if its $V(X)$-covolume is a sum with finitely many terms.
 Finally, using
Theorem~\ref{t:Scovolumes}, we now normalise the Haar measure
$\mu$ on $G=\Aut(X)$ so that for all uniform lattices $\G \leq G$,
the covolume of $\G$ is \[ \mu(\G \bs G) = \Vol(\G \quot V(X))\]

There are several constructions of uniform lattices acting on
polyhedral complexes.  For buildings, there are arithmetic
lattices, while for the Davis--Moussong complex, the Coxeter group
$W_L$ associated to $L$ is a uniform lattice
in $\Aut(X_L)$.  If a polyhedral complex $X$ is constructed as the
universal cover of a (faithful) finite complex of finite groups,
then the fundamental group of this complex of groups is a uniform
lattice acting without inversions on $X$.

\section{Covolumes for finitely many link types}\label{s:covolumes}

In this section we state and prove Theorem~\ref{t:covols}, a very general
restriction on covolumes of uniform lattices acting on complexes with finite
fundamental domains, and thus finitely many isometry classes of links.  Then in
Section~\ref{ss:examples} we give some applications of this theorem.

\subsection{Covolume restrictions}\label{ss:mainthm}

Theorem~\ref{t:covols} below can be seen as a generalisation of
the following result on tree lattices, due to Levich~(in
\cite{ros1:tctl}, Lemma 9.1.1).  Let $T_m$ be the $m$-regular
tree, and suppose $a/b$, a rational in lowest terms, is the
covolume of a uniform lattice $\G \leq \Aut(T_m)$. Then $b$ is not
divisible by any primes greater than $m$, and if $m$ is prime then
$b$ is not divisible by $m^2$ (such a $b$ is called an
$m$-number). The key step in the proof is showing that for all
vertices $v$ of the tree, the order of the stabiliser $\G_v$ is an
$m$-number. For this, $\G_v$ is injected into the automorphism
group of a finite rooted tree, and this rooted automorphism group
is computed. In the proof of Theorem~\ref{t:covols} we also inject
the vertex stabiliser into a rooted automorphism group, but derive
information about the order of the rooted group without actually
computing the group itself.

\begin{theorem}\label{t:covols} Let $X$ be a $d$-dimensional polyhedral complex
such that $G \bs X$ is finite, where $G = \Aut(X)$.  Let $L_1$,
$L_2$,\ldots, $L_m$ be the finite simplicial complexes such that
for each $x \in V(X)$, $\Lk(x,X)$ is isometric to some $L_i$.  Let
\[ p_1^{\alpha_1} p_2^{\alpha_2} \cdots p_M^{\alpha_M}\] be the
lowest common multiple of the orders of the groups $\Aut(L_i)$, $1
\leq i \leq m$, with each $p_j$ prime, $p_1 < p_2 < \cdots < p_M$,
and $\alpha_j > 0$ for $1 \leq j \leq M$.  For $1 \leq i \leq m$
and each $(d-1)$-cell $\sigma$ of $L_i$, let $\Fix_i(\sigma)$ be
the subgroup of $\Aut(L_i)$ which fixes that cell pointwise. Let
\[ p_1^{\alpha_1'}p_2^{\alpha_2'}\cdots p_M^{\alpha_M'}\] be the
lowest common multiple of the orders of the groups
$\Fix_i(\sigma)$, for $1 \leq i \leq m$, so that $0 \leq \alpha_j'
\leq \alpha_j$ for $1 \leq j \leq M$.  Then if $\Gamma$ is a
uniform lattice in $G$, its covolume $\mu(\G \bs G) = a/b$ is a rational in
lowest terms, such that $b$ is:

(a) not divisible by any primes other than $p_1, p_2, \ldots, p_M$; and

(b) if for some $j$ we have $\alpha_j' = 0$, then not divisible by
$p_j^{\alpha_j + 1}$. \end{theorem}

\begin{proof} By definition the covolume of $\Gamma$ is the finite sum \[
\mu(\G \bs G) = \sum_{x \in \G \bs V(X)} \frac{1}{|\Gamma_x|} \]
We will show that (a) and (b) hold for the order of each
stabiliser $\G_x$, and so complete the proof.

Since $\Gamma_x$ is finite, for $n$ sufficiently large $\Gamma_x$
injects into the finite group $H_n := \Aut(B(x,n))$. We will prove
by induction on $n\geq 1$ that (a) and (b) hold for the order of
$H_n$. By Lagrange's Theorem, (a) and (b) then hold for the order
of any subgroup of $H_n$.

To begin the induction, when $n = 1$, the ball $B(x,n)$ consists only of those
$d$-cells corresponding to $(d-1)$-cells in $\Lk(x,X)$. Hence, we may identify $H_1$ with
a subgroup of $\Aut(L_i)$ for some $i$. Thus (a) and (b) hold for the
order of $H_1$ (irrespective of the values of the $\alpha_j'$).

Assume (a) and (b) hold for $|H_n|$. An element of
$H_{n+1}=\Aut(B(x,n+1))$ fixes $x$, and so restricts to an element
of $H_n$.  Let $\varphi: H_{n+1} \rightarrow H_n$ be this
restriction homomorphism, with kernel $K$ and image $I$, so that
\[|H_{n+1}| = |K||I|\]  By induction, since $I$ is a subgroup of
$H_n$, (a) and (b) hold for the order of $I$.  We will show that
(a) holds for $|K|$, and that if $\alpha_j' = 0$ then no power of
$p_j$ divides $|K|$.

Now, since $B(x,n)$ contains finitely many cells, finitely many
vertices of $X$ lie in the boundary of $B(x,n)$. Enumerate these
boundary vertices as $x_1,x_2,\ldots,x_N$, and consider the
restriction homomorphisms \[K \rightarrow \Aut(B(x_k,1)), \mbox{
for } 1 \leq k \leq N \] We may identify each $\Aut(B(x_k,1))$
with a subgroup of $\Aut_\ell(L_i)$, for some $i$. Since (a) holds
for the order of a subgroup of any $\Aut_\ell(L_i)$, we have that
(a) holds for the order of the image of this restriction map.
Moreover, if $K'$ is any subgroup of $K$, then (a) also holds for
the order of the image of the restriction homomorphisms \[K'
\rightarrow \Aut(B(x_k,1)), \mbox{ for } 1 \leq k \leq N\]

For each $k$, at least one $(d-1)$-cell in the link of the vertex $x_k$
corresponds to a $d$-cell of $B(x,n)$. So, as elements of $K$ fix $B(x,n)$
pointwise, at least one $(d-1)$-cell in the link of each $x_k$ is fixed
pointwise by $K$. If for some $j$ we have $\alpha_j' = 0$, then $p_j$ does not
divide the order of any of the subgroups $\Fix_\ell(e,L_i)$. So if $\alpha_j' =
0$, the order of the image of the restriction homomorphism $K \rightarrow
\Aut(B(x_k,1))$ is not divisible by $p_j$.  And if $K'$ is any subgroup of $K$,
then the order of the image of the restriction map $K' \rightarrow
\Aut(B(x_k,1))$ will also not be divisible by $p_j$.

Put $K_0 = K$, and for $k=1,2,\ldots,N$ define $K_k$ and $I_k$ to
be respectively the kernel and image of the restriction map
\[K_{k-1} \rightarrow \Aut(B(x_k,1))\]  Then \[ K_N \leq K_{N-1}
\leq \cdots \leq K_2 \leq K_1 \leq K_0 = K \] This implies that
for all $1 \leq k \leq N$, (a) holds for the order of $I_k$, and
if $\alpha_j' = 0$ then $|I_k|$ is not divisible by $p_j$. Since
\[ |K| = |K_N| |I_N||I_{N-1}|\cdots|I_2| |I_1| \] and $K_N$ is
trivial, we are done. \end{proof}

\subsection{Examples}\label{ss:examples}

We now apply Theorem~\ref{t:covols} to some examples in dimensions 2, 3 and 4. 
The notation is as in the statement of the theorem, so that in each case $\mu(\G
\bs G) = a/b$ is a rational in lowest terms.

\begin{enumerate} \item Let $L_i = K_{m_i,n_i}$ be a complete bipartite graph
with $m_i \geq n_i \geq 2$, for $1 \leq i \leq m$.  Let $S_n$ be the symmetric
group on $n$ letters.  Then if $m_i = n_i$, \[\Aut(L_i) = (S_{m_i} \times
S_{n_i}) \rtimes S_2\] and if $m_i > n_i$, \[\Aut(L_i) = S_{m_i} \times
S_{n_i}\]  Theorem~\ref{t:covols} implies that $b$ is not divisible by any
primes greater than $M = \max_i\{m_i\}$. The subgroup of $\Aut(L_i)$ which
fixes an edge pointwise is $S_{m_i-1} \times S_{n_i-1}$.  Hence, if $M$ is
prime and the graph $K_{M,M}$ is one of the links then $b$ is not divisible by
$M^3$, and if $M$ is prime and $\max_i\{n_i\}<M$ then $b$ is not divisible by
$M^2$. \item Let $m=1$ and $L_1=L$ be the \emph{Petersen graph} (see, for
example,~\cite{o1:tg} pp. 240--241).  The 10 vertices of this graph may be
identified with the set of transpositions in $S_5$, and two  vertices are
joined by an edge if those transpositions are disjoint.  Now \[\Aut(L) =S_5\]
so $b$ is not divisible by any primes greater than 5. The subgroup of $\Aut(L)$
which fixes an edge pointwise is $S_2 \times S_2$, so $b$ is not divisible by
$3^2$ or $5^2$. \item Let $m=1$ and $L_1=L$ be the flag complex of the
projective plane over a finite field $\F_{q}$. The group $\Aut_0(L)$ of
type-preserving automorphisms of $L$ has index 2 in $\Aut(L)$, and is
isomorphic to $P\G L_3(\F_q)$, the group of incidence-preserving bijections of
the projective plane over $\F_q$. If $q = p^n$, where $p$ is prime, then the
order of $P\G L_3(\F_q)$ is given by~(see, for example, \cite{h1:pg} Theorem
2.8): \[|P\G L_3(\F_q)| = n|PGL_3(\F_q)| = nq^3(q^3-1)(q^2-1)\] Hence, $b$ is
not divisible by any primes other than those dividing $n$, $q$, $(q^2 + q + 1)$
and $(q \pm 1)$. Since $\Aut_0(L)$ acts transitively on the $(q+1)(q^2 + q +
1)$ edges of $L$, the subgroup of $\Aut_0(L)$ which fixes an edge pointwise has
order $nq^3(q-1)^2$. Depending on the value of $q$, this may tell us more about
the prime divisors of $b$. \item More generally, let $\cG$ be a finite, rank 2
Chevalley group over a finite field $\F_{q}$. Let $m=1$ and $L_1=L$ be the
spherical building associated to the $BN$-pair of $\cG$. The group $\Aut_0(L)$
is an extension of $\cG$ by $\Aut(\F_{q})$ (\cite{t1:bst}, Corollary 5.9). 
Table 1 of~\cite{gp1:ih} gives the orders of the groups $\cG$ and the number of
edges of the corresponding graphs $L$. Since $\Aut_0(L)$ acts transitively on
the set of edges of $L$, the order of a subgroup of $\Aut_0(L)$ which fixes an
edge of $L$ pointwise may thus be found.  The previous example is the case
where $\cG$ is of type $A_2$. \item Similar arguments to Example (3) may be
used for the 3-dimensional complexes constructed by Haglund--Paulin
in~\cite{hp1:ca} (see Example \ref{e:hp}, Section~\ref{ss:poly}), where the
link is the flag complex of projective 3-space.\item Let $m=1$ and $L_1 =L$ be
the first barycentric subdivision of the 120-cell, the regular polytope in
$\R^4$ whose boundary consists of 120 dodecahedrons (see, for
example,~\cite{s1:s120}).  The Davis--Moussong complex associated to $L$ is  a
4-dimensional right-angled hyperbolic building (\cite{js1:hcg}, proof of Theorem 2).  The
automorphism group of the 120-cell has order $120^2$~\cite{bp1:nn}, so $b$ has
no prime divisors other than 2, 3 and 5, and since the subgroup fixing a dodecahedron
pointwise is trivial, $b$ must be a factor of $120^2$.\end{enumerate}

\section{Covolumes for polygonal complexes with one link type}\label{s:onelink}

In the remainder of this note we consider 2-dimensional complexes.  A polygonal
complex is said to be an \emph{$(r,L)$-complex} if all of its 2-cells
are isometric to regular $r$-gons, and the links of all of its
vertices are isometric to a graph $L$.  In this section we
establish restrictions stronger than those of
Theorem~\ref{t:covols} on covolumes of uniform lattices acting on
certain $(r,L)$-complexes. First, in Sections~\ref{ss:rLcomplexes}
and~\ref{ss:uniformrL} we recall results on the existence and
uniqueness of these $(r,L)$-complexes, and constructions of
uniform lattices acting on them. Section~\ref{ss:covolsrL} then
contains our restrictions on covolumes.

\subsection{Existence and uniqueness of $(r,L)$-complexes}\label{ss:rLcomplexes}

In general, there may be uncountably many pairwise non-isomorphic
$(r,L)$-\linebreak complexes (see \cite{bb1:pc} Theorem 1.6, \cite{h2:pg},
and \cite{gp1:ih} Theorem 3.6). In the following cases, however,
local data does uniquely determine a polygonal complex.

\begin{enumerate}
\item Let $L = K_{m,n}$ be a complete bipartite graph, for $m, n
\geq 2$, and let $r \geq 5$, with $r$ even if $m \not = n$.  Let
$P$ be a regular right-angled hyperbolic $r$-gon. Then there
exists a unique connected $(r,L)$-complex $X$ such that all
2-cells of $X$ are isometric to $P$. This is due to Bourdon
(\cite{b1:ih}, Proposition~2.2.1) and Haglund (\cite{h1:ih},
Theorem~2). \item Let $L$ be the spherical building associated to
a finite Chevalley group $\cG$ of rank 2 over a finite field
$\F_q$. Then $L$ is a generalised $m$-gon, for some $m \geq 3$.
Let $r \geq 5$ and let $P$ be a regular hyperbolic $r$-gon with
all vertex angles $\frac{\pi}{m}$. A connected $(r,L)$-complex $X$
with all 2-cells isometric to $P$ is called an
\emph{$(r,L)$-building}.  We say that $X$ is \emph{locally
reflexive} if along each edge of $X$, the subcomplex consisting of
the 2-cells meeting that edge possesses an automorphism of order 2
(see~\cite{h1:ih} for the exact definitions).  Then if $q$ is
prime, and $r\geq 6$ is even, there exists a unique locally
reflexive $(r,L)$-building $X$. This is due to Haglund
(\cite{h1:ih}, Theorem~2).
\end{enumerate}

\subsection{Uniform lattices acting on $(r,L)$-complexes}\label{ss:uniformrL}

For each of the examples in Section~\ref{ss:rLcomplexes}, where an
$(r,L)$-complex is specified by local data, we describe a uniform
lattice acting on that $(r,L)$-complex. The constructions (due to
Bourdon~\cite{b2:if} and Gaboriau--Paulin~\cite{gp1:ih}) were
originally complexes of groups. Here, we state the stabilisers of
faces, edges and vertices in the quotient.

\begin{enumerate}
\item Let $X$ be as in (1) of Section~\ref{ss:rLcomplexes}. Then
there is a uniform lattice $\G \leq \Aut(X)$ so that the quotient
$\G\bs X$ is the polygon $P$. The face stabiliser is trivial. If
$L$ is $K_{m,m}$ then the stabiliser of each edge is $\Z/m\Z$ and
of each vertex is $\Z/m\Z \times \Z/m\Z$. If $L$ is $K_{m,n}$ with
$m \not = n$, so that $P$ has an even number of sides, then the
edge stabilisers alternate between $\Z/m\Z$ and $\Z/n\Z$, and the
vertex stabilisers are $\Z/m\Z \times \Z/n\Z$. This is a
rephrasing of \cite{b2:if}, Example 1.5(a). \item   Let $X$ be as
in (2) of Section~\ref{ss:rLcomplexes}. Then there is a uniform
lattice $\G \leq \Aut(X)$ so that the quotient $\G\bs X$ is the
polygon $P$. The face stabiliser is the group $B$ of the $BN$-pair
of the Chevalley group $\cG$. The edge stabilisers alternate
between $P_1$ and $P_2$, where $P_1$ and $P_2$ are parabolic
subgroups of $\cG$. The vertex stabilisers are the group $\cG$.
This is a rephrasing of~\cite{gp1:ih}, Section 3.1.4, Example (A).
\end{enumerate}

\subsection{Covolume restrictions}\label{ss:covolsrL}

We conclude by establishing restrictions on covolumes of uniform
lattices acting on the complexes described in
Section~\ref{ss:rLcomplexes}. If $X$ is as in (1) or (2) of
Section~\ref{ss:rLcomplexes}, then $\Aut(X)$ contains a finite
index normal subgroup $\Aut_0(X)$, the group of type-preserving
automorphisms, which acts without inversions.  Thus, any uniform
lattice $\Gamma \leq \Aut(X)$ has a finite index subgroup $\Gamma
\cap \Aut_0(X)$ which acts without inversions.  Now, in
Corollaries~\ref{c:bipartite} and~\ref{c:building} below, the
given sets of rational numbers are closed under multiplication by
positive integers, so we need only consider uniform lattices which
act without inversions.  Note that the lattices described in
Section~\ref{ss:uniformrL} act without inversions, since they were
constructed using complexes of groups.

The first result we will need is an analogue of
Theorem~\ref{t:covols}. Let $F(X)$ be the set of faces, or
2-cells, of a polygonal complex $X$. Theorem~\ref{t:facevols}
below gives a restriction on $F(X)$-covolumes.

\begin{theorem}\label{t:facevols} With the notation of Theorem~\ref{t:covols},
suppose $\Gamma$ is a uniform lattice which acts without
inversions. Then the $F(X)$-covolume of $\G$ is rational $a/b$ (in
lowest terms), such that $b$ is not divisible by any primes other
than $p_1,p_2,\ldots,p_M$. Moreover, if for some $j$ we have
$\alpha_j'=0$, then $b$ is not divisible by $p_j$.
\end{theorem}

Note that no power of $p_j$ can be a factor of $b$, in contrast to
Theorem~\ref{t:covols}.

\begin{proof} The $F(X)$-covolume is the finite sum \[ \Vol(\G \quot F(X)) =
\sum_{\sigma \in \G \bs F(X)} \frac{1}{|\Gamma_\sigma|} \] We
claim that the order of each stabiliser $\G_\sigma$ is not
divisible by any primes other that $p_1,p_2,\ldots,p_M$, and that
if for some $j$ we have $\alpha_j' = 0$, then the order of
$\Gamma_\sigma$ is not divisible by $p_j$.  The proof of this
claim is similar to that of the claim about orders of vertex
stabilisers in Theorem~\ref{t:covols}, except that we begin the
induction with the group of automorphisms without inversions of
$B(\sigma,0)$, which is trivial, so its order is not divisible by
any $p_j$.
\end{proof}

We will also use the following consequence of Theorem~\ref{t:Scovolumes}.  A
similar result holds for tree lattices, and is used by Rosenberg to establish a
restriction on the covolumes of uniform lattices acting on biregular trees
in~\cite{ros1:tctl}, Theorem 9.2.1.

\begin{corollary}\label{c:vertexfacevols} Let $X$ be a locally finite polygonal
complex such that $G\bs X$ is finite, where $G=\Aut(X)$. Then
there is a constant $c(X)$, depending only on $X$, such that for
all uniform lattices $\Gamma$ which act without inversions, \[
\mu(\G \bs G) = c(X)\Vol(\G \quot F(X))
\]
\end{corollary}

Therefore, if both $\mu(\G \bs G)$ and $\Vol(\G \quot F(X))$ are
known for just one uniform lattice $\G$ acting without inversions,
the constant $c(X)$ may be computed. Using the examples of uniform
lattices in Section~\ref{ss:uniformrL}, together with
Theorem~\ref{t:facevols}, we obtain the following results.

\begin{corollary}\label{c:bipartite} Let $X$ be as in (1) of Section~\ref{ss:rLcomplexes}, and $G=\Aut(X)$. Then if
$\G \leq G$ is a uniform lattice,
\[ \mu(\G \bs G) = \frac{r}{mn}\left( \frac{a}{b} \right)
\] where $a/b$ is rational (in lowest terms), and the prime divisors of $b$ are strictly less than $m$.
\end{corollary}

In fact, we show in~\cite{at1:cb} that every rational number of
the form given in Corollary~\ref{c:bipartite} can be obtained as
the covolume of some uniform lattice in $G$.

\begin{corollary}\label{c:building}  Let $X$ be as in (2) of Section~\ref{ss:rLcomplexes}, and $G=\Aut(X)$.
Then if $\G \leq G$ is a uniform lattice,
\[ \mu(\G \bs G) = \frac{r}{[\cG:B]}\left( \frac{a}{b} \right) \]
where $B$ is from the $BN$-pair of $\cG$, and $a/b$ is rational
(in lowest terms), such that $b$ is not divisible by any primes
other than those dividing the order of a subgroup of $\Aut_0(L)$
which fixes an edge of $L$ pointwise.
\end{corollary}

The value of $[\cG:B]$ can be computed from Table 1 of
~\cite{gp1:ih}.

\end{document}